\documentclass[12pt,english]{article}
\usepackage[T1]{fontenc}
\usepackage[latin9]{inputenc}
\usepackage[a4paper]{geometry}
\usepackage{amsmath}
\usepackage{amssymb}
\usepackage{graphicx}
\usepackage[authoryear]{natbib}
\usepackage{xargs}[2008/03/08]

\makeatletter
\usepackage{url}

\makeatother

\usepackage{babel}
\usepackage{listings}
\begin{document}
\title{UltimateKalman: Flexible Kalman Filtering and Smoothing Using Orthogonal
Transformations}
\author{Sivan Toledo\\
Blavatnik School of Computer Science, Tel Aviv University}
\maketitle
\begin{abstract}
UltimateKalman is a flexible linear Kalman filter and smoother implemented
in three popular programming languages: MATLAB, C, and Java. UltimateKalman
is a slight simplification and slight generalization of an elegant
Kalman filter and smoother that was proposed in 1977 by Paige and
Saunders. Their algorithm appears to be numerically superior and more
flexible than other Kalman filters and smoothers, but curiously has
never been implemented or used before. UltimateKalman is flexible:
it can easily handle time-dependent problems, problems with state
vectors whose dimensions vary from step to step, problems with varying
number of observations in different steps (or no observations at all
in some steps), and problems in which the expectation of the initial
state is unknown. The programming interface of UltimateKalman is
broken into simple building blocks that can be used to construct filters,
single or multi-step predictors, multi-step or whole-track smoothers,
and combinations. The paper describes the algorithm and its implementation
as well as a test suite of examples and tests.
\end{abstract}

\global\long\def\linearestimator{F}%

\global\long\def\modelfun{M}%

\global\long\def\penaltyfun{\phi}%

\global\long\def\objectivefun{\phi}%

\global\long\def\grad{\nabla}%

\newcommandx\hessian[1][usedefault, addprefix=\global, 1=]{\nabla_{#1}^{2}}%

\global\long\def\jacobian{\mathrm{J}}%

\global\long\def\exact#1{#1}%

\global\long\def\estimate#1{\hat{#1}}%

\global\long\def\controlpoint{\rho}%

\global\long\def\location{\ell}%

\global\long\def\noise{\epsilon}%

\global\long\def\observations{b}%

\global\long\def\expectation{\operatorname{E}}%

\global\long\def\iu{\mathbf{i}}%

\global\long\def\vecone{\mathbf{1}}%

\global\long\def\veczero{\mathbf{0}}%

\global\long\def\covold{\text{cov}}%

\global\long\def\nonjsscov{\operatorname{cov}}%

\global\long\def\cov{\operatorname{cov}}%

\global\long\def\var{\text{var}}%

\global\long\def\fim{\mathcal{I}}%

\global\long\def\loglikelihood{\mathcal{L}}%

\global\long\def\score{\mathcal{S}}%

\global\long\def\duration{\vartheta}%

\global\long\def\attenuation{a}%
\global\long\def\cmplxatt{\alpha}%

\global\long\def\initialphase{\varphi}%

\global\long\def\satclockerr{\eta}%

\global\long\def\ionodelay{\psi}%

\global\long\def\tropodelay{\xi}%

\global\long\def\xcorr{\operatorname{xcorr}}%

\global\long\def\diag{\operatorname{diag}}%

\global\long\def\rank{\operatorname{rank}}%

\global\long\def\erf{\operatorname{erf}}%

\global\long\def\erfc{\operatorname{erfc}}%

\global\long\def\range{\operatorname{range}}%

\global\long\def\trace{\operatorname{trace}}%

\global\long\def\ops{\operatorname{ops}}%

\global\long\def\prob{\operatorname{Prob}}%

\global\long\def\real{\text{Re}}%

\global\long\def\imag{\text{Im}}%

\global\long\def\square{\text{\ensuremath{\blacksquare}}}%

\global\long\def\irange{\boldsymbol{:}}%

\section{Introduction}

The invention of the Kalman filter by Rudolf  E. K{\'a}lm{\'a}n
in 1960~\citep{Kalman:1960:KalmanFilter} is considered one of the
major inventions of the 20th century. The filter efficiently and incrementally
tracks the hidden state of a linear discrete dynamic system; each
state estimate uses all the observations of the system up to that
point in time. The filter can also predict future states and with
suitable adaptations, can handle nonlinear dynamic systems and smooth
entire state trajectories. The literature on Kalman filters and their
applications is vast and the importance of Kalman filters is beyond
doubt. We mention a few relatively recent and relatively comprehensive
sources ~\citep{BrownHwangKalman3,GrewalAndrewsKalman4,AFreshLook2012},
but there are numerous other authoritative sources on Kalman filtering. 

Twelve years later, Duncan and Horn discovered that mathematically,
the Kalman filter computes the solution to a generalized linear least
squares problem~\citep{DucanHornKalman}. Algorithmically and numerically,
however, the Kalman filter algorithm is far from state-of-the art
algorithms for linear least squares problems. In particular, it does
not use orthogonal factorizations, which were invented and published
in the 1950s~\citep{Givens1954,HouseholderQR}. 

Numerically, the most stable algorithms for least squares minimization
are based on orthogonal transformations, and more specifically on
the QR factorization, the singular-value decomposition, or their variants~\citep{Bjorck:1996:NML,Golub:2013:MC,Higham:2002:ASN}.
Kalman filter algorithms are not, and many of them, including K{\'a}lm{\'a}n's
original algorithm, are based on algebraic building blocks, like explicit
matrix inversion, that are prone to instability.

In 1977, Paige and Saunders discovered and published an elegant Kalman
filter algorithm based on orthogonal transformations~\citep{PaigeSaunders:1977:Kalman}.
Their algorithm is a specialized QR factorization and it can easily
implement filtering, prediction, and smoothing.

Strangely, the Paige-Saunders algorithm appears to have had very limited
impact, even though it was about as efficient as other Kalman filters
and smoothers and was claimed to be more numerically stable. Their
paper was not cited much, nobody approached the authors to discuss
it~\citep{SaundersEmail2018}, and to the best of our knowledge,
it was never implemented (Paige and Saunders' paper describes the
algorithm and analyzes it, but does not mention an implementation).

The present paper and the software that it describes, called UltimateKalman,
aim to make the algorithm widely available and to highlight its advantages
over other Kalman filtering and smoothing algorithms. Indeed, the
Paige and Saunders algorithm is more flexible than other Kalman algorithms
in two important senses. First, unlike other Kalman filters, it does
not need to know the expectation of the initial state of the system.
Second, it can be easily generalized to handle quantities that are
added or dropped from the state vector. These features make modeling
easier, as we demonstrate with concrete examples in Sections~\ref{subsec:A-Projectile-Problem}
and~\ref{subsec:Clock-Offsets}. The algorithm can also easily handle
problems with a varying number of observations and with missing observations,
and it can both filter and smooth. Many other Kalman filter algorithms
lack these two characteristics, but some do posses them, so they are
not completely unique. 

UltimateKalman is not completely identical to the Paige-Saunders
algorithm, but rather a variant that is simpler and more general at
the same time. We explain later in the paper the differences from
the original algorithm, but at the same time acknowledge that all
the fundamental algorithmic ideas in UltimateKalman come from the
Paige-Saunders algorithm and paper.

Our implementation is split into a collection of easy-to-understand
building blocks from which a user can compose a variety of Kalman-based
computations, including filters, predictors, and smoothers, and combinations
of these. 

The implementation is available in three popular programming languages:
MATLAB, Java, and C. 

The rest of the paper is organized as follows. Section~\ref{sec:Background}
provides background material on discrete linear dynamic systems and
Kalman filters and smoothers. Section~\ref{sec:UltimateKalman} describes
the details of UltimateKalman. Section~\ref{sec:Implementation}
presents our implementations of the algorithm, and Section~\ref{sec:Examples-and-Tests}
describes a suite of examples and tests that come with UltimateKalman,
demonstrate its correctness, and show how to use it in various cases,
some nontrivial. We discuss the algorithm and the software and our
conclusions from this project in Section~\ref{sec:Discussion}.

This paper does not directly compare UltimateKalman to other Kalman
filtering algorithms and does not prove its correctness and its numerical
properties; the paper by Paige and Saunders addresses these issues
thoroughly, and the analyses there are equally applicable to UltimateKalman~\citep{PaigeSaunders:1977:Kalman}.

\section{Background}

\label{sec:Background}The discrete Kalman filter is a method to efficiently
estimate the state of a discrete linear dynamic system from indirect
observations. UltimateKalman can handle more general cases than Kalman
filters, so we describe here the more general version.

\subsection{Discrete linear dynamic systems}

The instantaneous state of a discrete dynamic system at time $t_{i}$
is represented by an $n_{i}$-dimensional \emph{state vector} $\exact u_{i}\in\mathbb{R}^{n_{i}}$.
We assume that $\exact u_{i}$ satisfies a recurrence that we refer
to as an \emph{evolution equation} and possibly another equation that
we refer to as an \emph{observation equation}. Note that we do not
require all the states to have the same dimension, although the uniform-dimension
case is very common. The evolution equation has the form
\begin{equation}
H_{i}\exact u_{i}=F_{i}\exact u_{i-1}+c_{i}+\epsilon_{i}\;,\label{eq:kalman-evolution-linear}
\end{equation}
where $H_{i}\in\mathbb{R}^{\ell_{i}\times n_{i}}$ and $F_{i}\in\mathbb{R}^{\ell_{i}\times n_{i-1}}$
are known full-rank matrices, $c_{i}\in\mathbb{R}^{\ell_{i}}$ is
a known vector, called a \emph{control vector}, that represents external
forces acting on the system, and $\epsilon_{i}$ is an unknown noise
or error vector. The control vector is often assumed to be the product
of a known matrix and a known vector, but this is irrelevant for the
Kalman filter. The noise or error vector $\epsilon_{i}$ admits state
vectors that do not satisfy the equation $H_{i}\exact u_{i}=F_{i}\exact u_{i-1}+c_{i}$
exactly. The matrix $H_{i}$ is often assumed to be the identity matrix,
but we do not require this (and do not require it to be square). When
$H_{i}$ and $F_{i}$ are square and full-rank, $u_{i}$ is a function
of $u_{i-1}$ and $c_{i}$ up to the error term. Obviously, the first
state $u_{0}$ that we model is not defined by an evolution recurrence.

Some of the state vectors $u_{i}$ (but perhaps not all) also satisfy
an \emph{observation} \emph{equation} of the form
\begin{equation}
o_{i}=G_{i}\exact u_{i}+\delta_{i}\;,\label{eq:kalman-observation-linear}
\end{equation}
where $G_{i}\in\mathbb{R}^{m_{i}\times n_{i}}$ is a known full-rank
matrix, $o_{i}\in\mathbb{R}^{m_{i}}$ is a known vector of observations
(measurements), and $\delta_{i}$ represents unknown measurement errors
or noise. The dimension $m_{i}$ of the observation of $u_{i}$ can
vary; it can be smaller than $n_{i}$ (including zero, meaning that
there are no observations of $u_{i}$), equal to $n_{i}$, or greater
than $n_{i}$. 

We can write all of the evolution and observation equations up to
step $k$ as a single large block-matrix equation, 
\begin{equation}
\left[\begin{array}{c}
o_{0}\\
c_{1}\\
o_{1}\\
c_{2}\\
\vdots\\
\vdots\\
c_{k}\\
o_{k}
\end{array}\right]=\left[\begin{array}{cccccc}
G_{0}\\
-F_{1} & H_{1}\\
 & G_{1}\\
 & -F_{2} & H_{2}\\
 &  & \ddots & \ddots\\
 &  &  & \ddots & \ddots\\
 &  &  &  & -F_{k} & H_{k}\\
 &  &  &  &  & G_{k}
\end{array}\right]\left[\begin{array}{c}
\exact u_{0}\\
\exact u_{1}\\
\exact u_{2}\\
\vdots\\
\exact u_{k-1}\\
\exact u_{k}
\end{array}\right]+\left[\begin{array}{c}
\delta_{0}\\
\epsilon_{1}\\
\delta_{1}\\
\epsilon_{2}\\
\vdots\\
\vdots\\
\epsilon_{k}\\
\delta_{k}
\end{array}\right]\;.\label{eq:kalman-assembly}
\end{equation}
We denote this system by $\observations=Au+e$. The matrix $A$ and
the vector $b$ are known. The noise or error terms $e$ are not known,
but we assume that they are small. Our task is to estimate $u$ from
$A$ and $b$. 

If $e$ is random, has zero expectation $\expectation(e)=0$, and
has a known nonsingular covariance matrix $\cov(e)=\expectation(ee^{T})$
with $\cov(e)^{-1}=U^{T}U$, where $U$ is square then the solution
$\estimate u$ of the \emph{generalized least squares }problem
\begin{eqnarray}
\hat{u} & = & \arg\min_{u}\left\Vert U\left(Au-b\right)\right\Vert _{2}^{2}\label{eq:gls}\\
 & = & \left(A^{T}U^{T}UA\right)^{-1}A^{T}U^{T}Ub\nonumber 
\end{eqnarray}
is the so-called best linear unbiased estimator (BLUE) of $u$~\citep{Aitken:1934:generalized-BLUE}.
If we add the assumption that $e$ has a Gaussian (normal) distribution,
then the same minimizer is also the maximum-likelihood estimator of
$u$. (In this paper, capital letters and small letters denote unrelated
objects; in this paragraph, for example, $\expectation$ denotes the
expectation and $e$ denotes an error or noise vector.)

If the structure of the $G_{i}$s, $F_{i}$s, and $H_{i}$s guarantees
that the rank of $A$ always equals the number of columns, the system
is called \emph{observable}. This guarantees that~(\ref{eq:gls})
has a unique solution. We use this term also in a more concrete sense:
we say that a matrix or a block of a matrix is \emph{observable} if
its rank equals its column dimension. We also use the terms \emph{flat
}to describe a rectangular matrix or block with more columns than
rows (a flat matrix cannot be observable) and \emph{tall} to describe
a rectangular matrix or block with more rows than columns. 

Because our goal is to estimate the states of a \emph{dynamic} system,
we also denote the linear system~(\ref{eq:kalman-assembly}) by $\observations^{(k)}=A^{(k)}u^{(k)}+e^{(k)}$,
to provide a notation for the matrix $A$ and for the right-hand $b$
at a particular step $k$. Similarly, we denote $\cov(e^{(k)})^{-1}=(U^{(k)})^{T}U^{(k)}$.

\subsection{Kalman filters and smoothers}

Kalman filters and smoothers are a large family of efficient algorithms
for solving problem~(\ref{eq:gls}) for some of the \emph{$\estimate u_{k}$}s
when $\cov(e)$ is block diagonal with known blocks, 
\[
\cov(e)=\left[\begin{array}{cccccc}
C_{0}\\
 & K_{1}\\
 &  & C_{1}\\
 &  &  & \ddots\\
 &  &  &  & K_{k}\\
 &  &  &  &  & C_{k}
\end{array}\right]\;.
\]
That is, we assume that the matrices
\begin{eqnarray*}
K_{i} & = & \cov\left(\epsilon_{i}\right)=\expectation\left(\epsilon_{i}\epsilon_{i}^{T}\right)\\
C_{i} & = & \cov\left(\delta_{i}\right)=\expectation\left(\delta_{i}\delta_{i}^{T}\right)
\end{eqnarray*}
are known and that the offdiagonal blocks of $\cov(e)$ are all zero:
\begin{eqnarray*}
\expectation\left(\epsilon_{i}\delta_{j}^{T}\right) & = & 0\;\text{for all}\;i\;\text{and}\;j\\
\expectation\left(\epsilon_{i}\epsilon_{j}^{T}\right) & = & 0\;\text{for}\;i\neq j\\
\expectation\left(\delta_{i}\delta_{j}^{T}\right) & = & 0\;\text{for}\;i\neq j\;.
\end{eqnarray*}
Most Kalman filtering and smoothing algorithms make additional assumptions.
For example, many Kalman filters assume that $G_{0}=I$, which is
equivalent to assuming that the expectation of $u_{0}$ is known.,
UltimateKalman also makes an additional assumption: that all the
$C_{i}$s and all the $K_{i}$s are nonsingular.

The generalized least-squares solution of Equation~(\ref{eq:kalman-assembly})
estimates all the state vectors $u_{1},\ldots,u_{k}$ using all the
observations up to and including step $k$. We denote the vectors
that make up this solution by 
\begin{equation}
\left[\begin{array}{c}
\estimate u_{0|k}\\
\estimate u_{1|k}\\
\vdots\\
\estimate u_{k|k}
\end{array}\right]\;.\label{eq:u-estimate-k}
\end{equation}
The vector $\estimate u_{k|k}$ is called the \emph{filtered} estimate
at step $k$. This estimate uses all the available observations of
present and past states, but not of any future state. Vectors $\estimate u_{0|k},\ldots,\estimate u_{k-1|k}$
are step-$k$ \emph{smoothed }estimates; they use observations of
past, present, and future states. We can hope to compute filtered
estimates almost in real time, whereas smoothed estimates can only
be computed after a time lag. If we extend system~(\ref{eq:kalman-assembly})
with one or more block rows and columns that represent only evolution
equations, the new vector components of the solutions are \emph{predicted
}estimates. For example, if we add a block row and a block column
that contain $-F_{k+1}$, $H_{k+1}$, and $c_{k+1}$, but not $G_{k+1}$
and not $o_{k+1}$, the last vector in the solution, denoted $\estimate u_{k+1|k}$,
is a prediction of $u_{k+1}$ from the information we have up to step
$k$. We can obviously continue to predict into the future by adding
more evolution equations.

Kalman filters are efficient incremental algorithms that produce filtered
and predicted estimates. Given $H_{k}$, $F_{k}$, $G_{k}$, $c_{k}$,
$o_{k}$, $C_{k}$, $K_{k}$, and the compact data structure that
was used to estimate $\estimate u_{k-1|k-1}$, Kalman filters quickly
compute $\estimate u_{k|k}$ and update the data structure~\citep{BrownHwangKalman3,GrewalAndrewsKalman4,AFreshLook2012,Kalman:1960:KalmanFilter}.
The data structure itself is of size $\Theta((n_{k-1}+n_{k})^{2})$
and the number of operations required is $O((n_{k-1}+m_{k}+n_{k})^{3})$.

\section{UltimateKalman and its heritage}

\label{sec:UltimateKalman}This section describes the UltimateKalman
algorithm. The algorithm is a slight simplification of the algorithm
of Paige and Saunders~\citep{PaigeSaunders:1977:Kalman} in that
UltimateKalman uses block orthogonal transformations whereas the
algorithm of Paige and Saunders uses Givens rotations.  The algorithm
is also a generalization of the algorithm of Paige and Saunders in
that we allow the user to specify $H_{i}$ and we allow the dimension
of the state vector to change from step to step.

\subsection{A specialized QR factorization}

UltimateKalman computes estimates of the state vectors using a thin
QR factorization of the weighted matrix

\begin{equation}
U^{(k)}A^{(k)}=\left[\begin{array}{cccccc}
W_{0}G_{0}\\
-V_{1}F_{1} & V_{1}H_{1}\\
 & W_{1}G_{1}\\
 & -V_{2}F_{2} & V_{2}H_{2}\\
 &  & \ddots & \ddots\\
 &  &  & \ddots & \ddots\\
 &  &  &  & -V_{k}F_{k} & V_{k}H_{k}\\
 &  &  &  &  & W_{k}G_{k}
\end{array}\right]\;,\label{eq:weighted-coefficient-block-matrix}
\end{equation}
where $W_{i}^{T}W_{i}=C_{i}^{-1}$ and $V_{i}^{T}V_{i}=K_{i}^{-1}$.
The factorization is computed using a series of orthonormal transformations
that are applied to block rows to reduce $U^{(k)}A^{(k)}$ to a block
upper triangular form: 
\[
R^{(k)}=\left(Q^{(k)}\right)^{T}\left(U^{(k)}A^{(k)}\right)\;,
\]
where
\begin{equation}
R^{(k)}=\left[\begin{array}{cccccc}
R_{0,0} & R_{0,1}\\
 & R_{1,1} & R_{1,2}\\
 &  & R_{2,2} & R_{2,3}\\
 &  &  & \ddots & \ddots\\
 &  &  &  & R_{k-1,k-1} & R_{k-1,k}\\
 &  &  &  &  & \tilde{R}_{k,k}
\end{array}\right]\;.\label{eq:block-R-factor}
\end{equation}
The diagonal blocks $R_{i,i}$ are normally square and upper triangular,
but are also allowed to be rectangular with more columns than rows.
The superdiagonal blocks $R_{i-1,i}$ are not identically zero. The
same series of transformations is applied to the weighted right-hand
side vector
\begin{equation}
U^{(k)}b^{(k)}=\left[\begin{array}{c}
W_{0}o_{0}\\
V_{1}c_{1}\\
W_{1}o_{1}\\
V_{2}c_{2}\\
\vdots\\
\vdots\\
V_{k}c_{k}\\
W_{k}o_{k}
\end{array}\right]\;.\label{eq:block-rhs}
\end{equation}
We denote the transformed right-hand side by $y^{(k)}=\left(Q^{(k)}\right)^{T}\left(U^{(k)}b^{(k)}\right)$.

The transformations are discarded immediately after they are applied
to the matrix and vector; no representation of $Q^{(k)}$ is stored.

\subsection{Observability}

Many Kalman algorithms rely on the assumption that $A$ is observable,
which guarantees that all the diagonal blocks $R_{i,i}$, as well
as $\tilde{R}_{k,k}$, are square and upper triangular. For example,
assuming that $G_{0}$ is square or tall (or is the identity) and
that the $F_{i}$s and $H_{i}$s are square, along with the standard
assumption that all of them have full rank, guarantees that $A$ is
observable. When all the diagonal blocks are square and triangular,
we can compute the estimates using back substitution, starting from
$\estimate u_{k|k}$ and ending with $\estimate u_{0|k}$.

However, UltimateKalman works and can provide useful estimates even
when $A$ is not always observable, and even when it is never observable.
We explain the different cases in terms of the structure of the $R$
factor and their meaning to the user.

If $\tilde{R}_{k,k}$ is flat but all the other $R_{i,i}$ blocks
are square and triangular, the system is \emph{not yet} observable,
but it may become observable in a future state. We currently do not
have enough observations to estimate any of the states. As the system
evolves and additional observations are made, the system may become
observable, allowing us to estimate the states that are currently
unobservable. 

If for some $i<k$ the diagonal block $R_{i,i}$ is flat, then states
$u_{0},\ldots,u_{i}$ are not observable and will never become observable.
For any assignment of states $u_{i+1}$ and up, there is a nontrivial
space of equally good (in the sense of~(\ref{eq:gls})) estimates
for $u_{0},\ldots,u_{i}$. UltimateKalman tolerates this situation
because even in this case, the observations of $u_{0},\ldots,u_{i}$
do provide useful information on future states.

If a state is not observable (or not yet observable), the method in
UltimateKalman that returns an estimate of the state vector return
a vector consisting of \texttt{NaN} values (not-a-number, a floating
point value that indicates that a value is not available). This informs
client code that the state is not observable.

\subsection{The Paige-Saunders factorization algorithm}

UltimateKalman uses the technique of the Paige-Saunders algorithm
to produce $R^{(k)}$ incrementally. 

Step $k$ starts by adding to $R^{(k-1)}$ a block row and a block
column that express an evolution equation, 
\[
\left[\begin{array}{cccc}
\ddots & \ddots\\
 & R_{k-2,k-2} & R_{k-2,k-1}\\
 &  & \tilde{R}_{k-1,k-1}\\
 &  & -V_{k}F_{k} & V_{k}H_{k}
\end{array}\right]\;.
\]
We now examine the block 
\begin{equation}
\left[\begin{array}{c}
\tilde{R}_{k-1,k-1}\\
-V_{k}F_{k}
\end{array}\right]\;.\label{eq:vconcat-Rtilde-VF}
\end{equation}
If this block is flat (or more generally, if its rank is smaller than
$n_{k}$, which implies that it can be orthogonally reduced to a flat
block), the algorithm leaves the bottom $2$-by-$2$ block as is,
but treats it as a $1$-by-$2$ block matrix $\begin{bmatrix}R_{k-1,k-1} & R_{k-1,k-1}\end{bmatrix}$,
with the same column partitioning. This allows the algorithm to handle
all shapes of~(\ref{eq:vconcat-Rtilde-VF}) in a uniform way. That
is, we rename the blocks as
\[
\left[\begin{array}{cc}
R_{k-1,k-1} & R_{k-1,k-1}\end{array}\right]=\left[\begin{array}{cc}
\tilde{R}_{k-1,k-1}\\
-V_{k}F_{k} & V_{k}H_{k}
\end{array}\right]\;.
\]
Otherwise, the algorithm computes a QR factorization of~(\ref{eq:vconcat-Rtilde-VF}),
uses the resulting $R$ factor as $R_{k-1,k-1}$, and applies the
orthonormal transformation to the last block column and to the right-hand
side $y$. This transforms the $R$ factor into
\[
\left[\begin{array}{cccc}
\ddots & \ddots\\
 & R_{k-2,k-2} & R_{k-2,k-1}\\
 &  & R_{k-1,k-1} & R_{k-1,k}\\
 &  &  & \bar{R}_{k,k}
\end{array}\right]\;.
\]

Block row $k-1$ is now \emph{sealed}; it will not change any more.
If $\tilde{R}_{k-1,k-1}$ was square, then so is $R_{k-1,k-1}$. If
$\tilde{R}_{k-1,k-1}$ was flat, then $R_{k-1,k-1}$ might be either
square or flat, depending on $V_{k}F_{k}$.

The bottom right block $\bar{R}_{k,k}$ is not upper triangular, and
it might be completely missing if~(\ref{eq:vconcat-Rtilde-VF}) is
square or flat. If $\bar{R}_{k,k}$ is square or tall and we now need
to predict $\estimate u_{k|k-1}$ (and perhaps additional future states),
we compute the QR factorization of $\bar{R}_{k,k}$ and apply the
transformation to $y$. We denote the $R$ factor of $\bar{R}_{k,k}$
by $\check{R}_{k,k}$. 

If there are no observations of $u_{k}$, then $\tilde{R}_{k,k}=\check{R}_{k,k}$
(or $\tilde{R}_{k,k}=\bar{R}_{k,k}$ if $\bar{R}_{k,k}$ is flat or
missing) and we are done with step $k$. 

If there are observations in step $k$, we add another block row to
the $R$ factor:
\[
\left[\begin{array}{cccc}
\ddots & \ddots\\
 & R_{k-2,k-2} & R_{k-2,k-1}\\
 &  & R_{k-1,k-1} & R_{k-1,k}\\
 &  &  & \bar{R}_{k,k}\\
 &  &  & W_{k}G_{k}
\end{array}\right]\;.
\]
If we computed $\check{R}_{k,k}$, in principle we can use it instead
of $\bar{R}_{k,k}$, but there is usually no significant benefit to
this. If the block

\[
\left[\begin{array}{c}
\bar{R}_{k,k}\\
W_{k}G_{k}
\end{array}\right]
\]
is flat, it becomes $\tilde{R}_{k,k}$ and we are done. Otherwise
we compute the QR factorization of this block, use the $R$ factor
as $\tilde{R}_{k,k}$, and apply the transformation to $y$. 

We now have $R^{(k)}$ and $y^{(k)}$. Note that we have not used
in this step block rows $1,\ldots,k-2$ of $R^{(k-1)}$ and $y^{(k-1)}$.

\subsection{\label{sec:forget-rollback}Forgetting and rolling back}

Given $R^{(k)}$ and $y^{(k)}$, the estimates $\estimate u_{i|k}$
are computed by back substitution, from the bottom up.

Therefore, if we need only filtered estimates $\estimate u_{k|k}$,
we only need the last sealed block row, row $k-1$, and the incomplete
block rows of step $k$. If we also need smoothed estimates but only
for time step $i+1$ and higher, we need to store block rows $i+1$
and higher, but not earlier rows.

Dropping old rows that would not be used for smoothing in the future
saves memory. UltimateKalman allows the user to \emph{forget} the
rows of steps $\leq i$ from $R^{(k)}$ and $y^{(k)}$. 

Smoothing uses $R^{(k)}$ but does not modify it. UltimateKalman
also retains $y^{(k)}$ when smoothing. This allows the algorithm
to smooth again later if more observations are obtained, enabling
easy implementation of strategies such as \emph{fixed-lag smoothing}~\citep{GrewalAndrewsKalman4},
in which each state is estimated once from past observations and from
observations of the next $n$ steps for some fixed lag $n$.

In many cases it is useful to predict future states before any observations
of them are available. In particular, by comparing the expectation
of the observations of a predicted state $G_{k}\estimate u_{k|k-1}$
with the actual observation vector $o_{k}$ it is sometimes possible
to detect and discard outlier observations.

UltimateKalman allows the user to predict future states while retaining
the ability to provide observations later. To do so, UltimateKalman
stores with sealed rows the incomplete diagonal block $\bar{R}_{k,k}$
and the associated right-hand side $\bar{y}_{i}$. If the user later
asks UltimateKalman to \emph{roll back} to step $k$, the algorithm
discards from memory sealed rows $k$ and higher and restores $\bar{R}_{k,k}$
and $\bar{y}_{i}$ as the bottommost incomplete row. Obviously, it
is not possible to roll back to a forgotten step.

\subsection{Computing the covariance matrices of estimates}

UltimateKalman computes representations of the covariance matrices
$\cov(\estimate u_{i|k})$ of estimates using orthogonal transformations
of $R^{(k)}$. We assume here that $\tilde{R}_{k,k}$ is square and
triangular (otherwise existing steps are not yet observable).

The covariance matrix of the filtered estimate satisfies
\[
\cov(\estimate u_{k|k})^{-1}=\tilde{R}_{k,k}^{T}\tilde{R}_{k,k}\;,
\]
so UltimateKalman simply returns $\tilde{R}_{k,k}$ as a representation
of $\cov(\estimate u_{k|k})$. We refer to this as an \emph{inverse-factor
}representation.

Producing inverse-factor representations of smoothed estimates requires
a series of orthogonal transformations. The algorithm first computes
the QR factorization of the bottom-right $2$-by-$1$ block of 
\[
R^{(k)}=\left[\begin{array}{cccc}
\ddots & \ddots\\
 & R_{k-2,k-2} & R_{k-2,k-1}\\
 &  & R_{k-1,k-1} & R_{k-1,k}\\
 &  &  & \tilde{R}_{k,k}
\end{array}\right]
\]
and applies the transformation to the entire two bottom block rows,
to produce
\[
\left[\begin{array}{cccc}
\ddots & \ddots\\
 & R_{k-2,k-2} & R_{k-2,k-1}\\
 &  & S_{k-1,k-1} & S_{k-1,k}\\
 &  & S_{k,k-1} & 0
\end{array}\right]
\]
with $S_{k,k-1}$ square. This block satisfies 
\[
\cov(\estimate u_{k-1|k})^{-1}=S_{k,k-1}^{T}S_{k,k-1}\;,
\]
so it is returned as a representation of the covariance matrix. The
algorithm now permutes the last two block rows to obtain
\[
\left[\begin{array}{cccc}
\ddots & \ddots\\
 & R_{k-2,k-2} & R_{k-2,k-1}\\
 &  & S_{k,k-1} & 0\\
 &  & S_{k-1,k-1} & S_{k-1,k}
\end{array}\right]
\]
and continues in the same way. The process continues with the QR factorization
of
\[
\left[\begin{array}{c}
R_{k-2,k-1}\\
S_{k,k-1}
\end{array}\right]\;,
\]
repeating the same procedure. The correctness of this algorithm is
shown in~\citep{PaigeSaunders:1977:Kalman} and in a perhaps slightly
clearer way, in~\citep{LocationEstimationFromTheGroundUp}. 

The computation proceeds upwards in $R^{(k)}$ and can produce the
covariance matrices of all the steps that have not been forgotten.

\section{Implementation}

\label{sec:Implementation}UltimateKalman is currently available
in MATLAB, C, and Java. This section describes the main features
of the implementations, but not their details. A user guide that is
distributed together with the source code of the library explains
the programming interfaces of the library and how to build it, how
to build test programs, and how to run these tests. 

Each implementation is separate and does not rely on the others. The
implementation includes MATLAB adapter classes that allow invocation
of the C and Java implementations from MATLAB and from GNU Octave,
a free open-source MATLAB-like environment. This allows a single
set of test functions to test all three implementations. The code
also includes standalone demonstration programs in Java and C, to
show users how to use the library without MATLAB.

The MATLAB implementation does not rely on any MATLAB toolbox, only
on functionality that is part of the core product. The implementation
also works under GNU Octave. The C implementation relies on basic
matrix and vector operations from the BLAS~\citep{BLAS,BLAS3ALG}
and on the QR and Cholesky factorizations from LAPACK~\citep{LAPACK-UG}.
The Java implementation uses the Apache Commons Math library for
both basic matrix-vector operations and for the QR and Cholesky factorizations. 

\subsection{Programming interface}

All three versions implement a data type that encapsulates the state
of a Kalman filter or smoother. The operations on the data type are
\texttt{evolve}, which describes an evolution equation (Equation~(\ref{eq:kalman-evolution-linear})
together with the associated covariance matrix $K_{i}$), \texttt{observe},
which describes an observation equation (Equation~(\ref{eq:kalman-observation-linear})
and $C_{i}$), \texttt{estimate}, which returns an estimated state
vector $\estimate u_{i}$ and its covariance matrix, \texttt{smooth},
which performs back substitution to compute smoothed estimates, and
\texttt{forget} and \texttt{rollback}, described in Section~\ref{sec:forget-rollback}.

The MATLAB and Java implementations use overloading (using the same
method name more than once, with different argument lists) to express
variants of these operations, while the C implementation uses \texttt{NULL}
to express missing or default values. The \texttt{estimate} method
returns two values in the MATLAB implementation, the estimated vector
and its covariance matrix, but only the estimated vector in the Java
and C implementations. The covariance matrix is returned in these
implementations by a separate method, \texttt{covariance}. These design
decisions reflect good programming practices in each programming language. 

\subsection{The representation of vectors and matrices}

The MATLAB implementation uses native MATLAB matrices and vectors.
The Java implementation uses the types \texttt{RealMatrix} and \texttt{RealVector}
from the Apache Commons Math library~(both are interface types with
multiple implementations). 

The C implementation defines a type called \texttt{matrix\_t} to
represent matrices and vectors. The implementation defines functions
that implement basic operations of matrices and vectors of this type.
The type is implemented using a structure that contains a pointer
to an array of double-precision elements, which are stored columnwise
as in the BLAS and LAPACK, and integers that describe the number
of rows and columns in the matrix and the stride along rows (the so-called
leading dimension in the BLAS and LAPACK interfaces). To avoid name-space
pollution, in client code this type is called \texttt{kalman\_matrix\_t}.

State vectors are not always observable. This topic is explained in
Section ~3.2 in the companion article. This situation usually arises
when there are not enough observations to estimate the state. The
function calls and methods that return estimates of state vectors
and the covariance matrices of the estimates return in such cases
a vector of \texttt{NaN}s (not-a-number, a floating point value that
indicates that the value is not available) and a diagonal matrix whose
diagonal elements are \texttt{NaN}.

\subsection{The representation of covariance matrices}

Like all Kalman filters, UltimateKalman consumes covariance matrices
that describe the distribution of the error terms and produces covariance
matrices that describe the uncertainty in the state estimates $\estimate u_{i}$.
The input covariance matrices are not used explicitly; instead, the
inverse factor $W$ of a covariance matrix $C=(W^{T}W)^{-1}$ is multiplied,
not necessarily explicitly, by matrices or by a vector. 

The programming interface of UltimateKalman expects input covariance
matrices $C_{i}$ and $K_{i}$ to be represented as objects belonging
to a type with a method \texttt{weigh} that multiplies the factor
$W_{i}$ such that $C_{i}=(W_{i}^{T}W_{i})^{-1}$ by a matrix $A$
or a vector $v$. In the MATLAB and Java implementations, this type
is called \texttt{CovarianceMatrix}. The constructors of these classes
accept many representations of a covariance matrix:
\begin{itemize}
\item An explicit covariance matrix $C$; the constructor computes an upper
triangular Cholesky factor $U$ of $C=U^{T}U$ and implements \texttt{X=C.weigh(A)}
by solving $UX=A$.
\item An inverse factor $W$ such that $W^{T}W=C^{-1}$; this factor is
stored and multiplied by the argument of \texttt{weigh}.
\item An inverse covariance matrix $C^{-1}$; the constructor computes its
Cholesky factorization and stores the lower-triangular factor as $W$.
\item A diagonal covariance matrix represented by a vector $w$ such that
$W=\text{diag}(w)$ (the elements of $w$ are inverses of standard
deviations).
\item A few other, less important, variants.
\end{itemize}
UltimateKalman always returns the covariance matrix of an estimated
state vector $\estimate u_{i}$ as an upper-triangular inverse factor
$W_{i}$. The MATLAB and Java implementations return covariance
matrices as objects of the \texttt{CovarianceMatrix} type (always
with an inverse-factor representation); the C implementation simply
returns the inverse factor as a matrix.

\section{Examples and tests}

\label{sec:Examples-and-Tests}We implemented an extensive set of
tests for UltimateKalman. The individual tests are implemented by
several MATLAB functions. Most of the functions receive as a first
argument a handle to a function that serves as a factory of UltimateKalman
filters. These functions are invoked by a top-level function called
\texttt{replication} that defines the factory function and performs
the tests. The factory function can produce objects of either the
MATLAB implementation or objects of adapter classes that invoke the
C or Java implementation.

The function that performs performance testing is a little different:
it includes a factory function and it can test multiple implementations,
to enable plotting their performance on one graph.

The tests generate graphs similar to the ones presented below. The
user can inspect the output visually, to ensure that the results are
similar to those presented here. The tests do not produce pass/fail
flags. 

The tests were run on a laptop with an Intel quad-core 14nm i7-8565U
processor running Windows~11 using MATLAB version R2021b. We also
verified that the MATLAB implementation works correctly under GNU Octave.
The C version is compiled into a MATLAB-callable dynamic link library
(a so-called \texttt{mex} file) by MATLAB itself using a script,
\texttt{compile.m}. In our tests, MATLAB used the C compiler from
Microsoft's Visual Studio 2019. The Java version is compiled using
Eclipse so that it can be used by Java~1.8 and up (this is the version
that MATLAB R2021b uses) and is packaged into a \texttt{jar} file
by a simple shell script, \texttt{build.bat}. 

\subsection{Basic tests}

\begin{figure}
\includegraphics[width=0.47\textwidth]{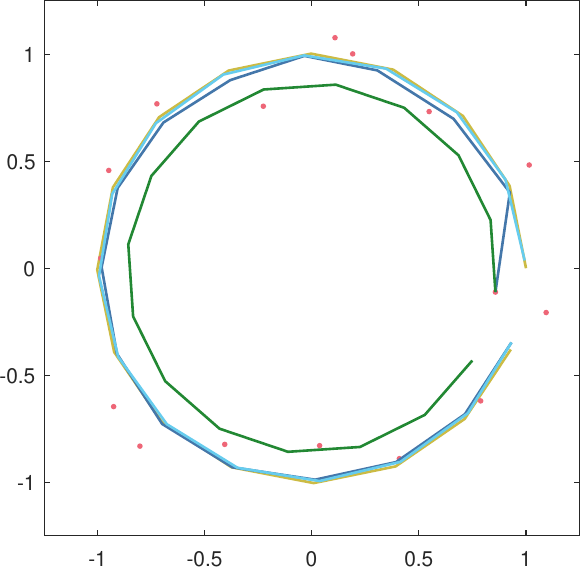}\hfill{}\includegraphics[width=0.47\textwidth]{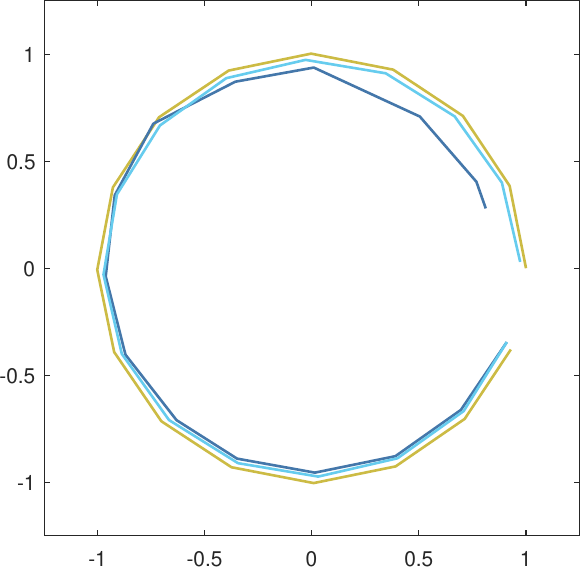}

\caption{\label{fig:rotation}Kalman filtering, prediction, and smoothing of
the trajectory of a rotating point in the plane. On the left $G=I$
and observations are indicated by red dots. The yellow curve represents
a simulated system; the blue curve represents filtered estimates,
the cyan whole-trajectory smoothed estimates, and the green predictions
from one observation. On the right $G=\begin{bmatrix}1 & 0\end{bmatrix}$;
we cannot predict the state from one such observation, and we cannot
produce a filtered estimate of the first state.}
\end{figure}
We demonstrate and test the basic features of UltimateKalman using
a simple model of a point in the plane that rotates around the origin.
The initial state is $\begin{bmatrix}1 & 0\end{bmatrix}^{T}$. The
evolution matrix is
\[
F=\left[\begin{array}{cc}
\cos(\alpha) & -\sin(\alpha)\\
\sin(\alpha) & \cos(\alpha)
\end{array}\right]
\]
for $\alpha=2\pi/16$, the observation matrix $G$ has between $1$
and $6$ rows, with an identity in the first two rows, $K=0.001^{2}I$
and $C=0.1^{2}I$. With these parameters, the rotation is very accurate,
but the observations are not. The system is simulated for $16$ steps,
just short of a complete rotation.

Figure~\ref{fig:rotation} shows the results of simulation, predicting,
filtering, and smoothing with $G=I$ and with $G=\begin{bmatrix}1 & 0\end{bmatrix}$.
The code, called \texttt{rotation}, first simulates the system and
produces the ground-truth $u_{1},\dots,u_{15}$ and the observations
$o_{0},\dots,o_{15}$. Then the code creates an UltimateKalman filter
and runs it for 16 steps while providing only the first observation
$o_{0}$. This attempts to predict $u_{1},\dots,u_{15}$. Then the
code rolls back to step~1 and runs the filter again, providing all
the observations. Finally, the code smooths the trajectory and collects
the smoothed observations.

The first observation in the case of $G=I$ is quite inaccurate, so
predictions from it are far from the real track, but they do follow
nicely the system dynamics of exact rotation. The filtered estimates,
on the other hand, improve quickly. The smooth estimates are nearly
perfect, at least visually.

When $G=\begin{bmatrix}1 & 0\end{bmatrix}$ the block $\bar{R}_{0,0}$
is flat, so the algorithm fails to produce a filtered estimate of
$u_{0}$, but it does produce filtered estimates of $u_{1},\dots,u_{15}$. 

This example is also used to test the code with overdetermined $G$s.

\subsection{Variance variations}

\begin{figure}
\includegraphics[width=0.47\textwidth]{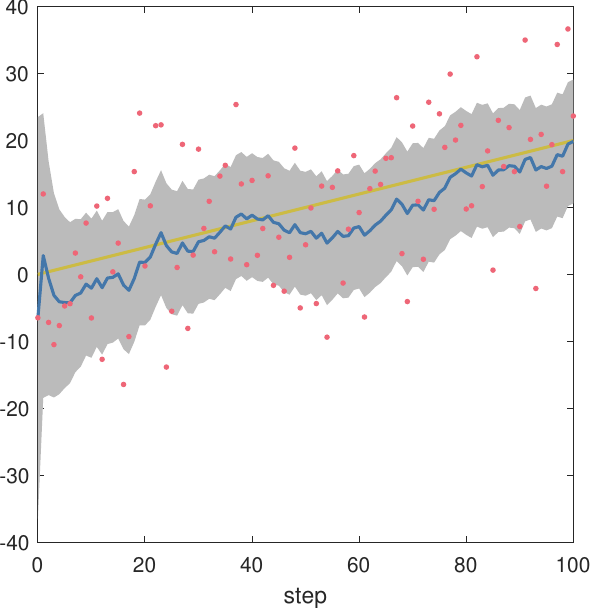}\hfill{}\includegraphics[width=0.47\textwidth]{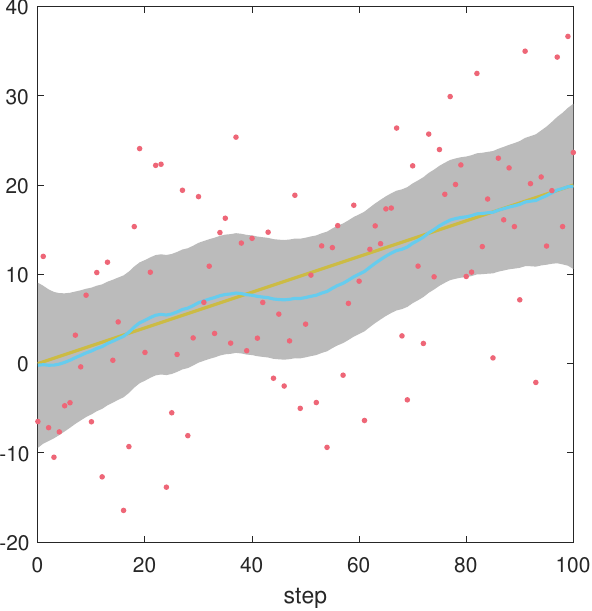}\\
\includegraphics[width=0.47\textwidth]{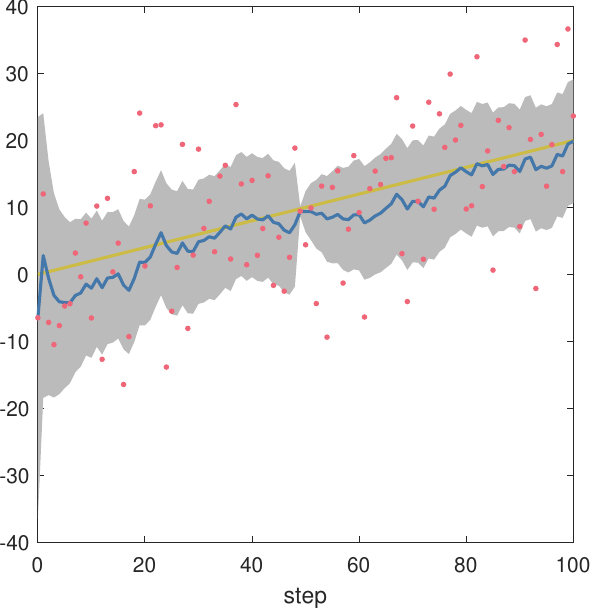}\hfill{}\includegraphics[width=0.47\textwidth]{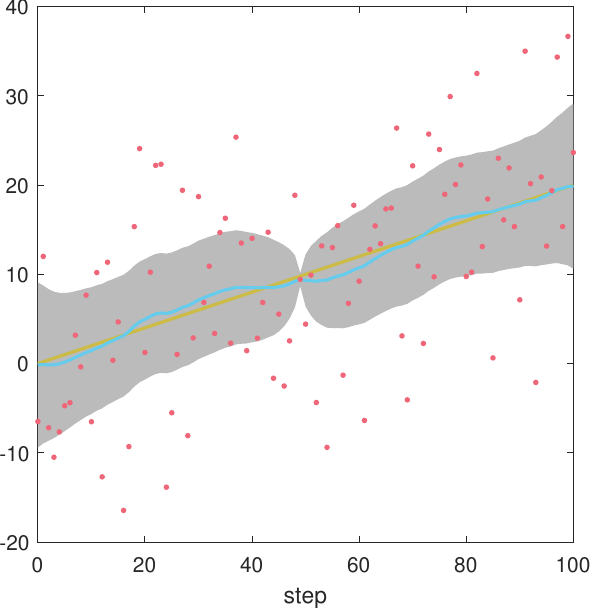}

\caption{\label{fig:constant-and-slope}Kalman filtering and smoothing of a
scalar trajectory. The yellow lines show the actual trajectory, the
red dots show the observations, the blue lines show the filtered estimates
(left column), the cyan lines show smoothed estimates (right column),
and the gray areas represent values that are $\pm3\protect\estimate{\sigma}_{i}$
from the estimate $\protect\estimate u_{i}$, where $\protect\estimate{\sigma}_{i}$
is the estimated standard deviation of $\protect\estimate u_{i}$.
In the top row, the standard deviation of $\delta_{i}$ is always
$10$; in the bottom row, it drops down to $0.25$ at iteration $50$
and then goes back up to $10$.}
\end{figure}
The next example demonstrates and tests the evaluation of covariance
matrices, and it also demonstrates the effectiveness of Kalman filtering
even when the model of the system is not perfect. 

The code  simulates a scalar that evolves either according to $u_{i}=u_{i-1}+\epsilon_{i}$
with $\epsilon_{1}\sim N(0,1)$ (the errors are distributed normally
with expectation $0$ and standard deviation $1$) or according to
$u_{i}=u_{i-1}+0.2$. However, even in the latter case, the Kalman
filter uses the $u_{i}=u_{i-1}+\epsilon_{i}$ evolution equations,
which do not reflect the true dynamics of the system. The observations
are direct, $o_{i}=u_{i}+\delta_{i}$, with $\delta_{i}$ having a
standard deviation of $10$ in almost all cases.

Figure~\ref{fig:constant-and-slope} shows results for the case $u_{i}=u_{i-1}+0.2$.
The graphs show both the estimates and their standard deviations.
We see that the variance of the filtered estimates drops quickly in
the first few steps and then stabilizes. If the standard deviation
of $o_{50}$ is much smaller than the rest, $0.25$, the variance
of the filtered estimate also drops at that step, but it climbs back
up. The variance of the smoothed estimates is impacted in both directions
around step $50$. It also increases towards iteration $0$ and $100$,
but not dramatically.

All the estimates track the trajectory nicely, even though the Kalman
filter and smoother use evolution equations that are different from
those of the actual dynamic system. This is an important reason that
Kalman filtering is so useful in practice: it often works well even
when it models the dynamic system only approximately.

The results for a simulation with $u_{i}=u_{i-1}+\epsilon_{i}$ are
similar and not shown in the paper.

\subsection{Adding and removing parameters}

\label{sec:add-remove}The \texttt{add\_remove} example shows how
to add and remove parameters from the state vector. In the first two
steps, the filter tracks the constant $1$ using the evolution equation
$u_{i}=u_{i-1}+\epsilon_{i}$ and observation equation $o_{i}=u_{i}+\delta_{i}$
with both error terms having a standard deviation of $0.1$. In the
third step, we add a dimension to the state; the first argument to
\texttt{evolve} is now $2$. This causes UltimateKalman to construct
and use a matrix 
\[
H_{i}=\left[\begin{array}{cc}
1 & 0\end{array}\right]\;,
\]
which represents a normal evolution of the first component of the
state while not using any historical information about the second,
new component. The evolution matrix $F_{2}$ remains a $1$-by-$1$
identity in this step, but $G_{2}$ is now a $2$-by-$2$ identity;
$o_{2}$ becomes $2$-dimensional as well. The actual value of the
second parameter is $2$ and the observations reflect that. In the
next iteration, $F_{3}$ grows to $2$-by-$2$. After two iterations
with a 2-dimensional state vector, we drop the first (original) component
of the state vector by calling \texttt{evolve} with a first argument
$1$ and by providing an explicit matrix 
\[
H_{4}=\left[\begin{array}{c}
0\\
1
\end{array}\right]\;.
\]
This matrix causes the filter to retain the second component of the
state and to drop the first. In this step $F_{4}$ remains $2$-by-$2$
but $G_{4}$ is $1$-by-$1$. In the next iteration $F_{5}$ shrinks
back to a $1$-by-$1$ identity. 

This demonstrates how to handle addition and removal of parameters
and tests that UltimateKalman handles these cases correctly. The
evolution and observation equations are very simple and track the
two parameters separately, but this is immaterial for the addition
and removal procedures.

\subsection{A Projectile problem}

\label{subsec:A-Projectile-Problem}
\begin{figure}[tb]
\includegraphics[width=0.47\textwidth]{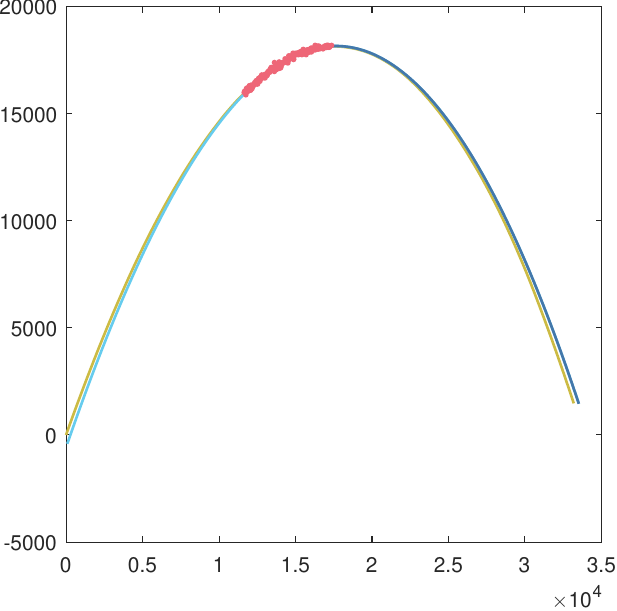}\hfill{}\includegraphics[width=0.47\textwidth]{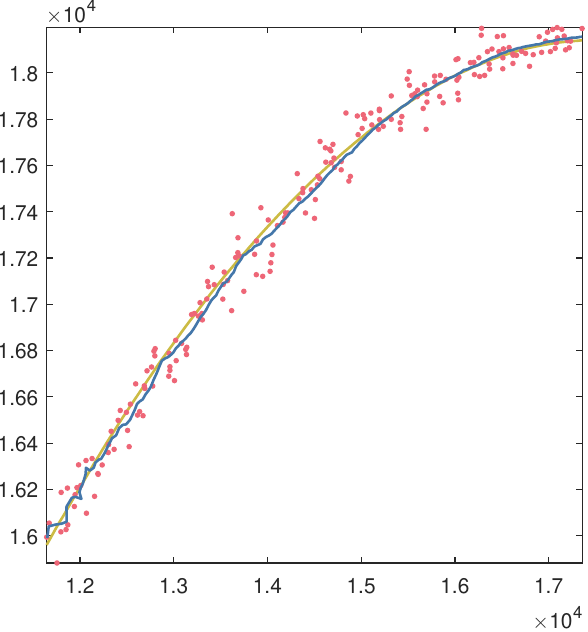}

\caption{\label{fig:projectile}Kalman filtering and smoothing of the trajectory
of a projectile affected by gravity and drag. Again, the yellow line
show the actual trajectory, the red dots show the observations, the
blue line shows the filtered estimates of the trajectory, and the
cyan line shows the smoothed estimates. Smoothing was performed on
the entire trajectory, to estimate the point of departure and the
point where the projectile hits the ground. The graph on the left
shows the entire trajectory and the graph on the right only the part
in which observations are available.}
\end{figure}
The next example that we consider is taken almost as is from Humpherys
et al.~\citep{AFreshLook2012}. The problem uses a linear dynamic
system to model a projectile. The states are four-dimensional; they
model the horizontal and vertical displacements and the horizontal
and vertical velocities. The model accounts for gravity, reducing
the vertical velocity by a constant every time step, and for drag,
which scales both components of the velocity by a constant in every
step. The matrices and vectors associated with the system are
\begin{eqnarray*}
F_{i}=F & = & \left[\begin{array}{cccc}
1 & 0 & \Delta t & 0\\
0 & 1 & 0 & \Delta t\\
0 & 0 & 1-b & 0\\
0 & 0 & 0 & 1-b
\end{array}\right]\quad\text{and}\quad c_{i}=c=\left[\begin{array}{c}
0\\
0\\
0\\
-9.8\Delta t
\end{array}\right]\\
G_{i}=G & = & \left[\begin{array}{cccc}
1 & 0 & 0 & 0\\
0 & 1 & 0 & 0
\end{array}\right]\;,
\end{eqnarray*}
where $\Delta t=0.1$ and $b=10^{-4}$. The example in~\citep{AFreshLook2012}
invites the reader to simulate the dynamic system for 1200 steps of
$0.1\,\text{s}$ each, starting from a known state, to generate noisy
observations of the displacements (not the velocities) in steps 400
to 600, and to estimate the trajectory using a Kalman filter. They
also ask the reader to use the filter to predict when and where the
projectile will fall back to its original altitude, and to estimate
the point of departure by reversing the dynamic system. 

We have successfully applied UltimateKalman to this problem. The
example code, \texttt{projectile}, generated the plots in Figure~\ref{fig:projectile},
which are similar to Figures~7.1a and~7.1b in~\citep{AFreshLook2012}.
(It seems that the measurement noise shown in their Figure~7.1b has
variance larger than the value 500 specified in~~\citep{AFreshLook2012}.) 

UltimateKalman also allows the user to easily extrapolate the trajectory
by evolving the filter without providing additional observations,
and to estimate the point of departure by smoothing time steps 0 to
600. Reversing the dynamic system is not required (and in particular,
it is not necessary to invert the evolution matrix). 

\subsection{Clock offsets in a distributed system}

\label{subsec:Clock-Offsets}
\begin{figure}[tb]
\hfill{}\includegraphics[width=0.47\textwidth]{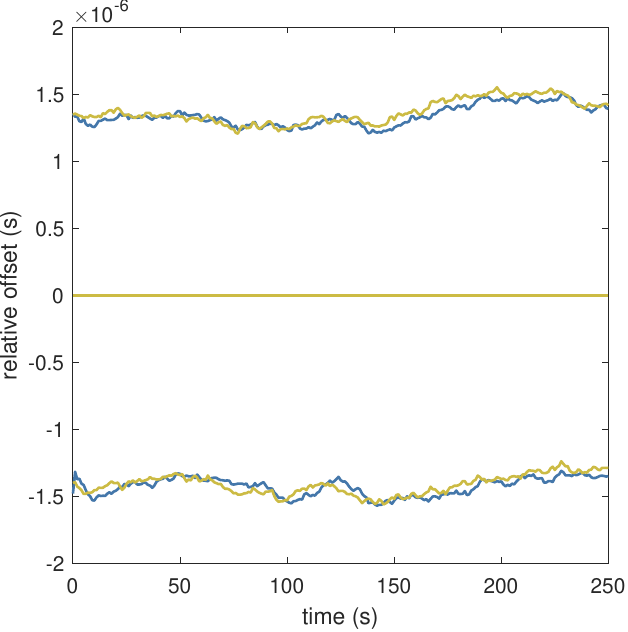}\hfill{}

\caption{\label{fig:clock-offsets}Kalman tracking of the relative offsets
of three clocks from times of arrival of radio packets. The graphs
show the relative offsets relative to the offset of the first clock.
The yellow lines are the simulated offsets and the blue lines are
the filtered estimates.}
\end{figure}
The \texttt{clock\_offsets} example highlights the utility of the
$H_{i}$ matrices. The aim is to estimate the relative offsets of
a set of clocks in a distributed system. Each clock is associated
with a receiver and we assume that at time $\tau$, clock $j$ shows
$t_{j}=\tau+f_{ij}$ where $f_{ij}$ is the offset of clock $j$ from
real time at the time in which it displays the value $t_{j}$. The
receivers receive radio packets from a beacon transmitter. The locations
of the transmitter and receivers are known, so the line-of-sight propagation
delays $d_{j}$ between the transmitter and receiver $j$ are also
known. The receivers estimate the time of arrival of the packets using
their local imperfect clocks. The observation equation for the time
of arrival of packet $i$ at receiver $j$ is
\[
t_{ij}=\tau_{i}+d_{j}+f_{ij}+\delta_{ij}\;,
\]
where $t_{ij}$ is the time-of-arrival estimate, as represented by
imperfect clock $j$, $\tau_{i}$ is the unknown time of departure
of the $i$th packet, $f_{ij}$ is the offset of clock $j$ at (local)
time $t_{ij}$, $d_{j}$ is the known delay to receiver $j$, and
$\delta_{ij}$ is the time-of-arrival estimation error.

The evolution equations are very simple: 
\[
f_{ij}=f_{i-1,j}+\epsilon_{ij}\;.
\]
They express the belief that the offsets change slowly. Note that
the number of evolution equations is equal to the number of clocks,
so is smaller than the dimension of the state vectors by one. The
code uses the following matrices and vectors, including an explicit
fixed $H_{i}$:
\begin{eqnarray*}
\left[\begin{array}{ccccc}
1 & 0 & 0 & 0 & 0\\
0 & 1 & 0 & 0 & 0\\
 &  & \ddots\\
0 & 0 & 0 & 1 & 0
\end{array}\right]\left[\begin{array}{c}
f_{i1}\\
f_{i2}\\
\vdots\\
f_{im}\\
\tau_{i}
\end{array}\right] & = & \left[\begin{array}{ccccc}
1 & 0 & 0 & 0 & 0\\
0 & 1 & 0 & 0 & 0\\
 &  & \ddots\\
0 & 0 & 0 & 1 & 0
\end{array}\right]\left[\begin{array}{c}
f_{i-1,1}\\
f_{i-1,2}\\
\vdots\\
f_{i-1,m}\\
\tau_{i-1}
\end{array}\right]+\left[\begin{array}{c}
\epsilon_{i1}\\
\epsilon_{i2}\\
\vdots\\
\epsilon_{1\ell}
\end{array}\right]\;,\\
\left[\begin{array}{c}
t_{i1}-d_{1}\\
t_{i2}-d_{2}\\
\vdots\\
t_{im}-d_{m}
\end{array}\right] & = & \left[\begin{array}{ccccc}
1 & 0 & \cdots & 0 & 1\\
0 & 1 & \cdots & 0 & 1\\
\vdots\\
0 & 0 & \cdots & 1 & 1
\end{array}\right]\left[\begin{array}{c}
f_{i-1,1}\\
f_{i-1,2}\\
\vdots\\
f_{i-1,m}\\
\tau_{i-1}
\end{array}\right]+\left[\begin{array}{c}
\delta_{i1}\\
\delta_{i2}\\
\vdots\\
\delta_{i\ell}
\end{array}\right]\;.
\end{eqnarray*}
The structure of $H_{i}$ and $F_{i}$ reflects the removal of $\tau_{i-1}$
and the introduction of $\tau_{i}$ in every step. 

The problem as presented up to now is clearly rank-deficient, because
the residual is invariant under an addition of a constant $T$ to
all the offsets and a subtraction of $T$ from all the depsarture
times. Nothing anchors the solution relative to absolute time; indeed,
we only want to estimate the relative offsets. To address this issue,
we add a pseudo-observation of one offset in the first step. This
removes the rank-deficiency. 

This model can be easily extended to allow for packets that are not
received by all receivers, for multiple transmitters, for slowly-changing
rate-errors instead of slowly-changing offsets, and so on. The model
can also allow receivers to join and leave the system without restarting
the estimation process, as we have done in Section~\ref{sec:add-remove}. 

The results are shown in Figure~\ref{fig:clock-offsets}.

\subsection{Performance testing}

\begin{figure}
\includegraphics[width=0.47\textwidth]{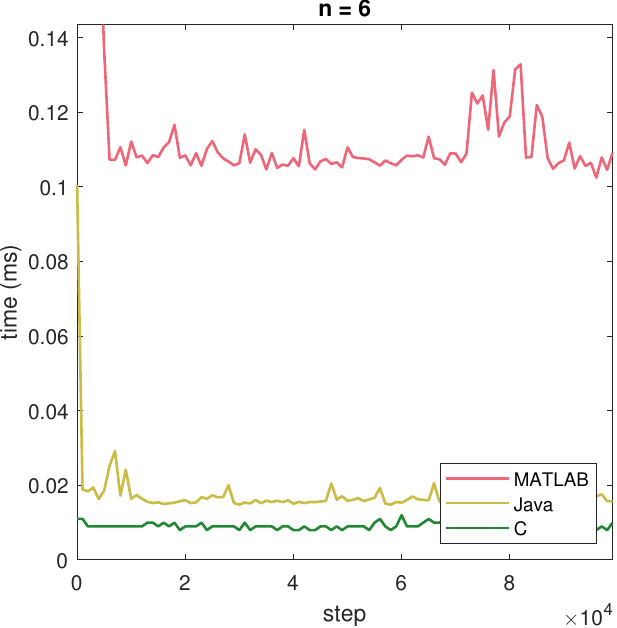}\hfill{}\includegraphics[width=0.47\textwidth]{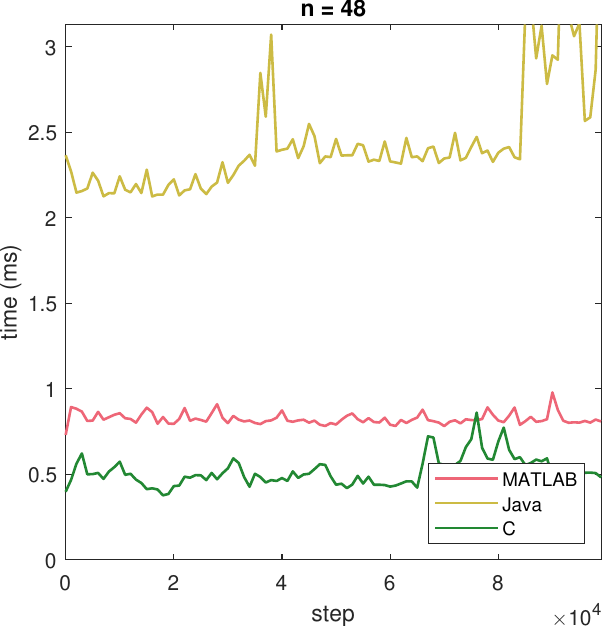}

\caption{\label{fig:stress-testing}Performance testing of the three UltimateKalman
implementations. The graphs show the time per step (window averages
of nonoverlapping groups of steps) when the algorithm filtered systems
with random square unitary matrices $F_{i}$ and $G_{i}$.}
\end{figure}

Figure~\ref{fig:stress-testing} shows the result of testing the
performance of the three implementations using their \texttt{perftest}
method. The experiments were carried out on a laptop with an Intel
i7 processor. The test uses two random square unitary matrices $F$
and $G$, sets $F_{i}=F$ and $G_{i}=G$, uses identities for $C_{i}$
and $K_{i}$, $c_{i}=0$, and a random Gaussian vector for $o_{i}$.
The matrices $H_{i}$ are identities created by the algorithm itself;
they are not passed as arguments. The use of unitary matrices avoids
overflows and underflows.

The graphs in Figure~\ref{fig:stress-testing} show the average runtimes
per step, averaged over groups of 1000 steps. 

The results show that on small problems ($n_{i}=6$), the C implementation
is the fastest, filtering taking about $8\,\mu\text{s}$ per step.
The Java implementation is about a factor of two slower, with periodic
jumps that are most likely caused by the garbage collector. The MATLAB
implementation is much slower, taking about $140\,\mu\text{s}$ per
step. Both the Java and MATLAB implementations are initially even
slower, probably due to just-in-time compilation of the code.

On larger problems ($n_{i}=48$), the C implementation is still the
fastest. The MATLAB implementation is now only about two times slower
than the C implementation. The improved ratio is most likely the
result of $\Theta(n^{3})$ dense-matrix operations (matrix multiplications
and QR factorizations) taking a significant fraction of the running
times. The data structure and method invocation overheads in C are
much smaller than in MATLAB, but dense-matrix operations are performed
at the same rate. On these larger problems, the Java implementation
is the slowest, most likely because of the relative poor performance
of the Apache Commons Math library relative to the BLAS and LAPACK
implementations that come with MATLAB.

Some of the changes in the runtimes that are visible in some of the
plots are likely due to the computer, a laptop, slowing down to avoid
overheating.

Smoothing takes more time per step, due to the back substitution phase,
and perhaps more importantly, a lot more memory, because the code
retains the blocks of $R$ and the transformed $b$. Still, the library
can smooth very long sequences. On the same computer we were able
to smooth 5 million steps with $n_{i}=m_{i}=6$ in about 96~s and
$100,000$ steps with $n_{i}=m_{i}=48$ in 84~s. Both runs required
less than 16GB of memory.

\section{Discussion}

\label{sec:Discussion}Orthogonal transformations are the bedrock
of numerical linear algebra. The numerical stability of algorithms
that rely solely or mostly on orthogonal transformations is often
both superior and easier to analyze than the stability of algorithms
that use non-orthogonal transformations or explicit matrix inversion.
From this standpoint, the fact that the Paige-Saunders algorithm did
not become the standard linear Kalman filter and smoother is both
puzzling and unfortunate. There is no good reason not to use it. 

UltimateKalman aims to rectify this defect. It is a simplified version
of the Paige-Saunders algorithm that is easier to implement, hopefully
also easier to understand, and is probably just as efficient. UltimateKalman
also generalizes the Paige-Saunders algorithm, making it more flexible.
While there are certainly many high-quality and well-documented Kalman
filter implementations, for example~\citep{KalmanFortran90,JSS:KalmanInR},
we are not aware of any other linear Kalman filter algorithm that
is as flexible and that can handle problems with varying state-vector
dimensions, does not require the expectation of the initial state,
can handle missing observations, and can easily filter, predict, and
smooth.

The Paige-Saunders algorithm does have two limitations that UltimateKalman
inherits. First, it assumes that the covariance matrices $C_{i}$
and $K_{i}$ are all nonsingular. Traditional Kalman filters can cope
with singular covariance matrices. However, this is not a significant
limitation, since physical observations are always associated with
some level of uncertainty, and physical evolution equations usually
are too. Continuous evolution equations can be exact, with no uncertainty
(e.g., equations of motion), but discrete evolution equations have
at least discretization (truncation) error terms. Even continuous
evolution equations sometimes suffer from uncertainty due to incomplete
physical modeling (e.g., a missing or simplified drag term). Second,
the algorithm is almost completely sequential, except for parallelism
within the matrix operations that are carried out in each step. In
real-time filtering applications, this limitation is meaningless,
as the application is inherently sequential. However, in post-processing
applications the ability to use many cores to reduce runtimes could
be beneficial; the Paige-Saunders algorithm and UltimateKalman lack
this ability.

The MATLAB implementation has been optimized for clarity and conciseness,
not for computational efficiency. In particular, the sequence of steps
in memory is represented by a dynamically-resized cell array. We made
this design choice in order to make UltimateKalman easy to port to
other languages and easy to specialize.

The C implementation was optimized for computational and memory efficiency
while retaining full generality (e.g., it can handle changes in the
dimension of the state vector, and it can smooth, like the other implementations)
and without a specialized memory manager. The main cost of these design
decisions is the dynamic allocation and deallocation of several matrices
in every step. Specializing the algorithm to simple special cases,
say only filtering and no changes in the state-vector dimension, would
have allowed the algorithm to avoid dynamic memory allocation, making
it faster and perhaps better suited to tiny embedded systems.

We hope that the numerical robustness afforded by the use of orthogonal
transformations, the flexibility, the convenient programming interfaces,
and the availability of the algorithm in multiple languages will make
UltimateKalman the standard linear Kalman filter and smoother in
new software. We also hope that authors and maintainers of statistical
software packages of wider scope~\citep{JSS:SpaceStateToolbooxMATLAB,JSS:SSpace},
which often include nonlinear Kalman filters and other state-space
models, will incorporate UltimateKalman into their packages, ideally
exposing some of the unique features through their own programming
and user interfaces. Finally, we hope that authors of code generators
that automatically generate Kalman filter code optimized for specific
cases~\citep{Autofilter2004} will also incorporate the algorithm
into their generators. The code clarity that we strived to achieve
should simplify such efforts, as well as efforts to implement the
algorithm in additional programming languages that are widely used
for numerical computing like Python and Julia.

\section*{Acknowledgements}

This research was supported in part by grant 1919/19 from the Israel
Science Foundation and by a grant from the Israeli Ministry of Science
and Technology..

\bibliographystyle{ACM-Reference-Format}
\bibliography{kalman}

\end{document}